\documentclass[12pt]{amsart}
\usepackage{graphicx}
\usepackage{amssymb}
\usepackage{epstopdf}
\theoremstyle{plain}

\newtheorem{theorem}[subsection]{Theorem}

\newtheorem{corollary}[subsection]{Corollary}
\newtheorem{conjecture}[subsection]{Conjecture}

\theoremstyle{definition}

\newtheorem{definition}[subsection]{Definition}

\newcommand\eps{\varepsilon}

\newcommand{\Z}{{\mathbf{Z}}}

\newcommand\N{{\mathbb{N}}}

\newcommand{\beq}{\begin{equation*}}
\newcommand{\eeq}{\end{equation*}}

\title{Roth type theorems in finite groups}
\author{J\'ozsef Solymosi}
 \address{University of British Columbia, Canada}
   \email{solymosi@math.ubc.ca}
\begin{document}

\begin{abstract}
We prove Roth type theorems in finite groups. Our main tool is the Triangle Removal Lemma of Ruzsa and Szemer\'edi.  
\end{abstract}

\dedicatory{Dedicated to the memory of Yahya Ould Hamidoune}

\thanks{This research was partially supported by an NSERC grant.}

\date{}

\maketitle

\section{Inroduction}
Brown and Buhler \cite{BB} have generalized Roth's classical theorem \cite{Ro} proving that every large subset of an abelian group contains three-term arithmetic progressions.  They proved the following.

\begin{theorem}\label{BBthm}
Let $V$ be an $n$-dimensional affine space over the field with $p^d$ elements, $p \neq 2$. Then for every $\eps > 0$ there is an $n(\eps)$ such that if $n = dim(V) \geq n(\eps)$ then any subset of $V$ with more than $\eps|V|$ elements must contain 3 collinear points (i.e., 3 points lying in a one-dimensional affine subspace).
\end{theorem} 

Frankl, Graham, and R\"odl gave a simple proof for Theorem \ref{BBthm} in \cite{FGR} using  a result of Ruzsa and Szemer\'edi, the so called (6,3)-theorem \cite{RSz}. This result is often called the Triangle Removal Lemma. In this paper we show that using the (6,3)-theorem and its extensions (see Theorem \ref{Removal} later) one can prove Roth-type theorems in non-abelian finite groups.  

The non-abelian case was considered by Bergelson, McCutcheon, and Zhang for amenable groups in \cite{BMZ}.
They proved that in any subset $E\subset G\times G$ which is of positive upper density (for the exact statement and definitions please refer to the article), one can find configurations of the form $\{(a, b), (ga, b), (ga, gb)\}.$  They also proved that if $G \times G \times G$ is partitioned into finitely many cells, one of these cells contains configurations of the form $\{(a, b, c), (ga, b, c), (ga, gb, c), (ga, gb, gc)\}.$ A stronger version of these results can be found in \cite{BM}.

Extending earlier work of Green \cite{Gr} on abelian groups, Kr\'al', Serra, and Vena \cite{KSV} proved Roth's theorem for finite groups. They also used the Triangle Removal Lemma as their main tool.

\begin{theorem}\label{Oriol}\cite{KSV}
Let $G$ be a finite group of odd order $N$ and $A$ a subset of its elements. If the number of solutions of the equation 
$xz = y^2$ with $x, y, z \in A$ is $o(N^2)$, then the size of $A$ is $o(N).$
\end{theorem}

\section{Results}
In this section we state some theorems and a conjecture which are similar in form to the results for amenable groups mentioned earlier. 
\begin{theorem}\label{elso}
For every $c>0$ there is a threshold $m_0\in \N$ such that if $H$ is a subgroup of a finite group $G$ and $|H|\geq m_0,$ 
then the following holds. Any set $S\subset G\times G$ with $|S|\geq c|G|^2$ contains three elements $(a,b), (ad,b),(a,db)$ where $d\in H.$ 
\end{theorem}
The proof of Theorem \ref{elso} uses very basic properties of finite groups only. It  was proved by Shkredov for abelian groups in \cite{Sh} where he gave an efficient bound on the density needed for such a configuration. 
He proved that there is a $\delta > 0$ such that, for any abelian group $G,$ if  $S\subset G\times G$ and $|S|\geq |G|^2/(\log\log|G|)^{\delta}$ then $A$ contains a triple $(a,b),(a+d,b),(a,b+d).$
Our method is not as effective since we are using the Triangle Removal Lemma for which the best known bound is still a $\log^*$ type bound \cite{Fo}.

One can use the  Triangle Removal Lemma in various ways to get different versions of Roth's Theorem in finite groups.

\begin{theorem}\label{harmadik}
For every $\delta>0$ there is a threshold $n_0\in \N$ such that if $G$ is a finite group of order $|G|\geq n_0$
then the following holds. Any set $S\subset G\times G$ with $|S|\geq \delta|G|^2$ contains three elements $(a,b), (a,c),(e,f)$ such that $ab=ec$ and $ac=ef.$ 
\end{theorem}

\noindent
It is possible that a stronger version of Theorem \ref{harmadik} holds, however I have not been able to prove it. 

\begin{conjecture}
For every $\delta>0$ there is a threshold $n_0\in \N$ such that if $G$ is a finite group of order $|G|\geq n_0$
then the following holds. Any set $S\subset G\times G$ with $|S|\geq \delta|G|^2$ contains four elements $(a,b),$ $(a,c),$ $(e,c),$ $(e,f)$ such that $ab=ec$ and $ac=ef.$ 
\end{conjecture}

\noindent
By choosing $d=e^{-1}a,$ the proof of Theorem \ref{harmadik} implies the following corollary.

\begin{corollary}\label{cor}
For every $\delta>0$ there is a threshold $m_0\in \N$ such that if $H$ is a subgroup of a finite group $G$ where $|H|\geq m_0,$ 
then the following holds. Any set $S\subset G\times G$ with $|S|\geq \delta|G|^2$ contains three elements $(a,b), (a,db),(ad^{-1},d^2b)$ where $d\in H.$ 
\end{corollary}
Note that the second coordinates of the three elements form an arithmetic progression, showing that any dense subset of $G$ contains a three term arithmetic progression, unless $d$ has order two or $e=a.$ The $e=a$ case would mean that $(a,b)=(a,c)=(e,f)$ so we can exclude it.  But one should consider the $d^2=1$ case as it is not avoidable in general. Elementary abelian 2-groups (Boolean groups) have no three-term arithmetic progressions since $b=d^2b$  there. That is why in Theorem \ref{BBthm} and \ref{Oriol} the order of the group is required to be odd. Theorem \ref{harmadik} (or rather its corollary) allows us to choose $d$ from a subgroup. So, if $G$ contains a sufficiently large subgroup without an order 2 element and we choose $d$ from there, then $b, db, d^2b$ form a non-degenerate arithmetic progression. So, Roth's theorem is true in a finite group if its order is divisible by a large odd number. Equivalently, Roth's theorem holds in finite groups where the index of the 2-Sylow subgroup  is large. The index of the $p$-Sylow subgroup of $G$ is denoted by $|G:Syl_p(G)|.$

\begin{theorem}[Roth's Theorem for Finite Groups]\label{RothGroups}
For every $c>0$ there is a threshold $m_0\in \N$ such that if $S$ is a subset of a finite group $G,$  $|S|\geq c|G|,$ and $G$ has a subgroup $H$ so that $|H:Syl_2(H)|\geq m_0$ then $S$ contains three distinct elements $b, db, d^2b$ where $d\in H.$ 
\end{theorem}

The next theorem illustrates the power of the Hypergraph Removal Lemma (Theorem \ref{Removal}). One important application of the Hypergraph Removal Lemma is to give a combinatorial proof \cite{So} of a multidimensional version of Szemer\'edi's Theorem proved by F\"urstenberg and Katz\-nel\-son \cite{FK}: Every dense subset of the $d$-dimensional integer grid, $\Z^d,$ contains a {\em corner}, i.e. $d+1$ points with coordinates 

$
\begin{array}{rlccrl}
(&a_1 , & a_2 , &\ldots  &, a_d&)\\
(&a_1+\delta , & a_2 &\ldots  &, a_d&)\\
(&a_1 , & a_2+\delta , &\ldots  &, a_d&)\\ 
\vdots\\
(&a_1 , & a_2 , &\ldots &, a_d+\delta&)
\end{array} 
$

\noindent
for some nonzero integer $\delta$. A similar statement holds for finite abelian groups. This is a generalization of Shkredov's Theorem \cite{Sh}. 

\begin{theorem}\label{negyedik}
For every $c>0$ there is a threshold $m_0\in \N$ such that if $H$ is a subgroup of a finite abelian group $G$ and $|H|\geq m_0$ 
then the following holds. Any set $S\subset \underbrace{G\times G\ldots \times G}_d$ with $|S|\geq c|G|^d$ contains $d+1$ elements with coordinates

$
\begin{array}{rlccrl}
(&a_1, & a_2, &\ldots  &, a_d &)\\
(&a_1+\delta, & a_2, &\ldots  &, a_d &)\\
(&a_1, & a_2+\delta, &\ldots  &, a_d &)\\ 
\vdots\\
(&a_1, & a_2, &\ldots &, a_d+\delta &)
\end{array} 
$

\noindent
where $\delta\in H.$ 
\end{theorem}

\section{Proofs}
In the proofs we will use the graph and the hypergraph removal lemmas. These results were proved by finding the right notation of regularity.
The idea of graph regularity played a key role in Szemer\'edi's proof for the Erd\H os-Tur\'an conjecture \cite{Sz}. The Regularity Lemma \cite{Sz2} became an important tool in graph theory and in additive combinatorics.  

\begin{definition}
An $r$-uniform hypergraph on $n$ vertices is defined by a subset of the $r$-element subsets (edges) of an $n$-element set (vertices). It is denoted by $\mathcal{H}^{r}_n$. If an $r$-uniform hypergraph on $r+1$ vertices contains all possible $r+1$ edges then it is called a {\em clique}.
\end{definition}

\begin{theorem}[Hypergraph Removal Lemma]\label{Removal}
If every edge of $\mathcal{H}^{r}_n$ is the edge of exactly one clique then it is sparse, i.e. its number of edges is $o(n^r).$
\end{theorem}

The $r=2$ case is the Triangle Removal Lemma \cite{RSz}. The $r=3$ case was proved by Frankl and R\"odl \cite{FR}. The general theorem was proved by R\"odl, Nagle, Shacht, and Skokan \cite{NRS,RS} and by Gowers \cite{G} independently. 

We are going to use the following quantitative version of the $r=2$ case. (Similar statements hold for larger $r$-s as well.)

\begin{theorem}[Triangle Removal Lemma]\label{triangle}
For every $\delta>0$ there is a $\delta'>0$ such that if a graph on $n$ vertices has $\delta n^2$ pairwise edge-disjoint triangles then it has at least $\delta'n^3$ triangles.
\end{theorem} 

\medskip

\noindent
{\em Proof of Theorem \ref{elso}}: We omit the condition that $d\in H,$ first. Let us define a tripartite graph on three vertex classes, $G_1,G_2,$ and $G_3$, where each $G_i$ is a copy of the group $G.$   In this tripartite graph three vertices $g_1\in G_1, g_2\in G_2$ and $g_3\in G_3$ span a triangle if $(g_1,g_2)\in S$ and $g_1g_2=g_3$. These are pairwise edge-disjoint triangles. The number of edge-disjoint triangles is  at least $\delta |G|^2$ so, by the Triangle Removal Lemma if the group is large enough\footnote{{\em large enough} means that $n$ is large enough to satisfy the right inequality in $\delta n^2\leq |S|<\delta'n^3$ where $\delta$ and $\delta'$ are from Theorem \ref{triangle} and $n=3|G|$.}  then there are more triangles which are spanned by the triangles defined by $S.$  One can find a triangle spanned by vertices (group elements) $a\in G_1,b\in G_2,$ and $c\in G_3$ so that $ab\not = c.$ As the three vertices span a triangle, one can find $(a,b),(x,b),(a,y)\in S$ so that $xb=c$ and $ay=c$. Then $(a,b), (cb^{-1},b), (a, a^{-1}c)\in S.$ Rewriting it we have $(a,b), (aa^{-1}cb^{-1},b), (a, a^{-1}cb^{-1}b)\in S.$ In order to complete the proof of Theorem \ref{elso} we have to show that $a^{-1}cb^{-1}$ can be chosen from $H.$ Let us consider the left and the right cosets of $H$. By the pigeon-hole argument there are elements $\ell, r\in G$ such that $|(\ell H \times Hr) \cap S| \geq \delta|H|^2.$ By choosing $m_0$ large enough (as earlier in the proof we selected $n$ being large enough to apply the Triangle Removal Lemma) one can repeat the previous argument for the tripartite graph on the three vertex classes $ \ell H, Hr$ and $\ell Hr.$ 
With this choice of vertex classes, $a\in \ell H, b\in Hr,$ and $c\in \ell Hr,$ so $a^{-1}cb^{-1}\in H.$ \qed

\medskip
\noindent
In the proof of Theorem \ref{harmadik} we will iterate the previous argument twice. In order to use the Triangle Removal Lemma twice, we need the quantitative version as stated in Theorem \ref{triangle}. A new element of the proof is that the subgroup $H$ is required to be abelian. A classical result of Erd\H os and Strauss \cite{ES} states that every finite group of order $n$ contains an abelian subgroup of order at least $\log{n}.$ This bound was improved significantly by Pyber \cite{P} who proved that every finite group contains an abelian subgroup of order at least $e^{c\sqrt{\log{n}}}$ for some $c.$ (This bound is best possible up to the constant $c.$)

\medskip

\noindent
{\em Proof of Theorem \ref{harmadik}}:  Let us choose a large abelian subgroup of $G$ denoted by $H.$ By the earlier mentioned theorems we can choose $H$ such that the order of $H$ is as large as necessary (depending on $\delta $).   Let us consider again  the left and the right cosets of $H$. There are elements $\ell, r\in G$ such that $|(\ell H \times Hr) \cap S| \geq \delta|H|^2.$ As before, we define a tripartite graph on the three vertex classes $ \ell H, Hr$ and $\ell Hr.$ In this tripartite graph three vertices $g_1\in \ell H, g_2\in Hr$ and $g_3\in \ell Hr$ span a triangle if $(g_1,g_2)\in S$ and $g_1g_2=g_3$. The number of edge-disjoint triangles is  at least $\delta|H|^2$ so by the Triangle Removal Lemma there are at least $\delta'|H|^3$ triangles in this tripartite graph. If $|H|$ is large enough then most of these triangles are spanned by vertices (group elements) $a_i\in\ell H,b_i\in Hr,$ and $c_i\in \ell Hr$ so that $a_ib_i\not = c_i.$ As the three vertices span a triangle, one can find $(a_i,b_i),(x_i,b_i),(a_i,y_i)\in S$ so that $x_ib_i=c_i$ and $a_iy_i=c_i$.  There are at least $c'|H|^3$ such triangles, therefore one can find an $x$ so that the number of triples $(a_i,b_i),(x,b_i),(a_i,y_i)\in S$ so that $xb_i=c_i$ and $a_iy_i=c_i$ is at least $\delta'|H|^2.$ The element $x$ and any pair of the triple $a_i,a_ib_i,c_i$ defines the triple $(a_i,b_i),(x,b_i),(a_i,y_i)$ uniquely. For the second part of the proof we define a new graph on the vertex set  $ \ell H \cup \ell Hr.$ The edges are defined by the pairwise edge-disjoint triangles spanned by $a_i,a_ib_i,c_i$ (for the fixed $x$).
For large enough $|H|$ we can apply the Triangle Removal Lemma again. It means that there is a triangle $T=(A,B,C)$ which is not given by $a_i,a_ib_i,c_i$ for some $i.$ $T$ is a triangle, so there are triples $(a_i,b_i),(x,b_i),(a_i,y_i),$ $(a_j,b_j),(x,b_j),(a_j,y_j),$ and $(a_k,b_k),(x,b_k),(a_k,y_k)$ triples defining edges $AB,$ $BC,$ and $AC$ (respectively) as follows; $AB=(a_i, a_ib_i),$ $BC=(a_jb_j,c_j),$ and $AC=(a_k,c_k).$  

Now we are ready to put together everything needed for the proof. By the selection of $T$ we have $a_i=a_k,$ $a_ib_i=a_jb_j,$ and $c_j=c_k.$ We show that if we choose the three elements  $(a,b)=(a_i, a_i^{-1}c_i),$ $(a,c)=(a_k,a_k^{-1}c_k),$ and $(e,f)=(a_j,a_j^{-1}c_j)$ from $S$ then these will satisfy the identities required. (Note that if $(a_s,c_s)$ was an edge in the first graph then we know that $(a_s, a_s^{-1}c_s)\in S.$) The second identity is $a_ja_j^{-1}c_j=a_ka_k^{-1}c_k$ which holds as $c_j=c_k.$ The first one requires to show that $a_ja_k^{-1}c_k=a_ia_i^{-1}c_i.$ Note that $c_i=xb_i,$  $c_k=c_j=xb_j,$ and $a_k=a_i,$ so one needs to show that $a_ja_i^{-1}xb_j=xb_i.$ Let us write $a_i=\ell \alpha_i,$ $a_j=\ell \alpha_j,$ $x=\ell \alpha_x,$ $b_i=\beta_ir,$ and $b_j=\beta_jr$ where $\alpha_i,\alpha_j,\alpha_x,\beta_i,\beta_j\in H.$ Now the equation is equivalent to $\alpha_j\alpha_i^{-1}\alpha_x\beta_j=\alpha_x\beta_i.$ As $H$ is an abelian subgroup, this reduces to $\alpha_i\beta_i=\alpha_j\beta_j$ which holds because  $a_ib_i=a_jb_j.$ 
\qed

\medskip
\noindent
{\em Proof of Corollary \ref{cor} and Theorem \ref{RothGroups}.} In the previous proof we had the freedom to choose the abelian subgroup $H$ as we want, we only used that it was large enough. If the order of the group is divisible by a large odd number $M$ then by the prime factorization of $M$ and Sylow's theorem we know that it has a $p$-subgroup, for some odd $p,$ of order at least $\log M.$ According to an old result of Miller \cite{B,M} every $p$-group $P$ contains an abelian subgroup of size at least $c\sqrt{|P|}$ . If we choose $H$ to be this subgroup then the three second coordinates $b,c,f$ (or $b, db, d^2b$ in Corollary \ref{cor} and in Theorem \ref{RothGroups}) are distinct and they form an arithmetic progression of length 3.    
\qed

\medskip

\noindent
{\em Proof of Theorem \ref{negyedik}}:  Let us find group elements $a_1, a_2, \ldots, a_d$ such that $|(a_1+H)\times (a_2+H)\times\ldots \times (a_d+H)\cap S|\geq c|H|^d.$ We define a $(d+1)$-partite $d$-uniform hypergraph. The vertex sets are $V_i=a_i+H$ for $1\leq i\leq d$ and $V_{d+1}=a_1+a_2+\ldots +a_d+H.$ Every $e\in \{(a_1+H)\times (a_2+H)\times\ldots \times (a_d+H)\cap S\}$ is an edge of the hypergraph. These are the {\em generator edges}, they will generate the hypergraph as follows. If a generator edge has coordinates $e=(v_1,v_2,\ldots ,v_d)$ then the $d+1$ vertices $v_1,v_2,\ldots ,v_d,$ and $ \sum_{i=1}^dv_i$ span a clique. These are pairwise edge-disjoint cliques since from the sum, $ \sum_{i=1}^dv_i,$ and all coordinates of $e$ but one, we can always find the missing coordinate uniquely. As there are $c|H|^d$ such cliques, by the Hypergraph Removal Lemma there should be another clique.  Let us denote the vertices of this clique by $w_1,w_2,\ldots ,w_{d+1}.$ We know that $\sum_{i=1}^dw_i\not=w_{d+1}.$ Since $w_1,w_2,\ldots ,w_{i-1},w_{i+1},\ldots ,w_{d+1}$ is an edge for any $1\leq i\leq d,$ then there is a generating edge $w_1,w_2,\ldots ,w_{i-1},x_i,w_{i+1},\ldots ,w_{d}$ such that $w_1+w_2+\ldots +w_{i-1}+x_i+w_{i+1}+\ldots +w_{d}=w_{d+1}.$ The $d$ generating edges $w_1,w_2,\ldots ,w_{i-1},x_i,w_{i+1},\ldots ,w_{d}$ and $w_1,w_2,\ldots ,w_{d}$ give the points of the corner as required. The value of $\delta$ is $\delta=w_{d+1}-\sum_{i=1}^dw_i$ which belongs to $H$ as the $a_i$-s cancel out. 
\qed

\section{Acknowledgements} The author is thankful to Mikl\'os Ab\'ert, P\'eter P\'al P\'alfy, and L\'aszl\'o Pyber for the useful conversations and to the referee for the careful reading and helpful suggestions. This research was supported by an NSERC grant.

\end{document}